\documentclass[11pt]{amsart}
\usepackage{amsfonts,amsmath} 
\usepackage[dvips]{graphicx} \newcommand{\halmos}{\rule{1ex}{1.4ex}}
\usepackage{psfrag}
\newcommand{\proofbox}{\hspace*{\fill}\mbox{$\halmos$}}
\newdimen\margin   
\def\COMMENT#1{}
\def\TASK#1{}

\def\noproof{{\unskip\nobreak\hfill\penalty50\hskip2em\hbox{}\nobreak\hfill%
       $\square$\parfillskip=0pt\finalhyphendemerits=0\par}\goodbreak}

\newcommand{\eps}{\varepsilon}

\newcommand{\Hy}{\mathcal{H}}

\newcommand{\Pa}{{\mathcal P}}
\newcommand{\G}{{\mathcal G}}

\newcommand{\n}{N_h}
\newcommand{\cn}{{\mathcal N}_h}
\newcommand{\F}{{\mathcal F}}
\newcommand{\N}{{\mathcal N}}
\newcommand{\B}{{\mathcal B}}
\newcommand{\Typ}{{\rm Typ}}
\newcommand{\Atyp}{{\rm Atyp}}
\newcommand{\Usef}{{\rm Usef}}

\newcommand{\C}{\mathcal{C}}

\newcommand{\R}{\mathcal{R}}
\newcommand{\D}{\mathcal{D}}

\newtheorem{firsttheorem}{Proposition}
\newtheorem{fact}[firsttheorem]{Fact}
\newtheorem{theorem}[firsttheorem]{Theorem}
\newtheorem{lemma}[firsttheorem]{Lemma}

\newtheorem{definition}[firsttheorem]{Definition}
\newtheorem{proposition}[firsttheorem]{Proposition}

\begin{document}
\title{$3$-uniform hypergraphs of bounded degree have linear Ramsey
numbers} \author{Oliver Cooley, Nikolaos Fountoulakis, Daniela K\"uhn and Deryk Osthus}

\begin{abstract} \noindent
Chv\'atal, R\"odl, Szemer\'edi and Trotter~\cite{CRST} proved that
the Ramsey numbers of graphs of bounded maximum degree are linear in
their order. We prove that the same holds for $3$-uniform
hypergraphs. The main new tool which we prove and use is an
embedding lemma for $3$-uniform hypergraphs of bounded maximum
degree into suitable $3$-uniform `pseudo-random' hypergraphs.
\end{abstract}

\maketitle 

\noindent
{\footnotesize
keywords: hypergraphs; regularity lemma; Ramsey numbers; embedding problems}

\vspace{0.4cm}

\section{Introduction}\label{intro}

\subsection{Ramsey numbers}
The \emph{Ramsey number $R(H)$} of a graph $H$ is defined to be the
smallest $N \in \mathbb{N}$ such that in every colouring of the
edges of the complete graph on $N$ vertices with two colours one can
find a monochromatic copy of $H$. In general, the best upper bound
on $R(H)$ is exponential in~$|H|$. However, if $H$ is sparse, then
one can sometimes improve considerably on this. A central result in
this area was proved by Chv\'atal, R\"odl, Szemer\'edi and
Trotter~\cite{CRST}. They showed that for every $\Delta$ there exists
a constant~$C=C(\Delta)$ such that all graphs~$H$ with maximum degree at
most $\Delta$ satisfy $R(H) \le C |H|$.

Here we prove an analogue of this result for $3$-uniform hypergraphs $\Hy$
of bounded maximum degree. Thus we now consider hyperedges (each
consisting of $3$ vertices) instead of edges. The \emph{degree of a vertex~$x$} 
in~$\Hy$ is defined to be the number of hyperedges which contain~$x$. The
\emph{maximum degree} $\Delta (\Hy)$ and the \emph{Ramsey number} $R(\Hy)$ of a
$3$-uniform hypergraph $\Hy$ are then defined in the obvious way.
\begin{theorem} \label{Ramseythm}
For every $\Delta\in\mathbb{N}$ there exists a constant~$C=C(\Delta)$ such that
all $3$-uniform hypergraphs~$\Hy$ of maximum degree at most~$\Delta$
satisfy $R(\Hy) \le C |\Hy|$.
\end{theorem}
\noindent
Kostochka and R\"odl~\cite{RK_OW} showed that Ramsey numbers of $k$-uniform
hypergraphs of bounded maximum degree are `almost linear' in their
orders. More precisely, they showed that for all $\varepsilon
,\Delta,k>0$ there is a constant $C$ such that $R(\Hy) \le
C|\Hy|^{1+\varepsilon}$ if $\Hy$ has maximum degree at most~$\Delta$.
Also, Haxell et al.~\cite{haxell1, haxell2} asymptotically
determined the Ramsey numbers of $3$-uniform tight and loose cycles.
In general, just as for graphs it follows immediately from Ramsey's
theorem that there is an absolute constant $C$ such that every
$3$-uniform hypergraph~$\Hy$ satisfies $R(\Hy) \le C^{|\Hy|}$.
Moreover, a well-known probabilistic argument due to Erd\H{o}s also
gives an exponential lower bound in $|\Hy|$ if $\Hy$ is complete.%

\subsection{Embedding graphs and hypergraphs}

The proof in~\cite{CRST} which shows that graphs of bounded maximum
degree have linear Ramsey numbers proceeds as follows: Given a
red/blue colouring of the edges of the complete graph on $N$
vertices, we consider the red subgraph $G$ and apply Szemer\'edi's
regularity lemma to it to obtain a vertex partition of $G$ into a
bounded number of clusters such that almost all of the bipartite
subgraphs induced by the clusters are `pseudo-random'. We now define a
reduced graph $R$ whose vertices are the clusters and any two of them
are connected by an edge if the corresponding bipartite subgraph of
$G$ is `pseudo-random'. Since $R$ is very dense, by Tur\'an's theorem
it contains a large clique~$K$. We now define an edge-colouring
of~$K$ by colouring an edge red if the density of the corresponding
bipartite subgraph of $G$ is large, and blue otherwise. An
application of Ramsey's theorem now gives a monochromatic clique of
order $k: =\Delta(H)+1$ in~$K$. Without loss of generality, assume
it is red. This corresponds to a large complete $k$-partite
subgraph~$G'$ of $G$ where all the bipartite subgraphs induced by the
vertex classes are `pseudo-random'. Since the chromatic number of
the desired graph~$H$ is at most~$k$, one can use the
`pseudo-randomness' of $G'$ to find a copy of $H$ in
$G'$. The tool which enables the final step is often called the
`embedding lemma' or `key lemma' (see e.g.~\cite{KSi}).

In our proof of Theorem~\ref{Ramseythm}, we adopt a similar
strategy. Instead of Szemer\'edi's regularity lemma, we will use the
regularity lemma for~$3$-uniform hypergraphs due to Frankl and
R\"odl~\cite{FR}. However, this has the problem that the
`pseudo-random' hypergraph into which we aim to embed our given
$3$-unifom hypergraph~$\Hy$ of bounded maximum degree could be very sparse
and not as `pseudo-random' as one would like it to be. This means that
the proof of the corresponding embedding lemma is considerably more
difficult and rather different to that of the graph version, while
the adaption of the other steps is comparatively easy. Thus we view the
embedding lemma (Lemma~\ref{embcor}) as the main result of this paper and also believe
that it will have other applications besides
Theorem~\ref{Ramseythm}. Its precise formulation needs some
preparation, so we defer its statement to
Section~\ref{sec:stateemb}.

\subsection{Organization of the paper}

In Section~\ref{sec:stateemb} we state the embedding lemma (Lemma~\ref{embcor}).
Our proof proceeds by induction on the order of the hypergraph we aim to embed.
This argument yields a significantly stronger result (Lemma~\ref{emblemma}).

In Section~\ref{sec:tools} we state several results which are important in the proof of the
induction step for Lemma~\ref{emblemma}. In particular, we will need a variant of the 
counting lemma,
which implies that for any $3$-uniform hypergraph~$\Hy$ of bounded size every suitable
`pseudo-random' hypergraph~$\G$ contains roughly as many copies of~$\Hy$ as one would expect in
a random hypergraph. We will also need an extension lemma, which states that for any
$3$-uniform hypergraph~$\Hy'$ of bounded size, any induced subhypergraph~$\Hy\subseteq \Hy'$ 
and any suitable
`pseudo-random' hypergraph~$\G$, almost all copies of~$\Hy$ in~$\G$ can be extended to
approximately the same number of copies of~$\Hy'$ as one would expect if~$\G$ were
a random hypergraph. In Section~\ref{uppercount} we derive our variant of the counting lemma
from that of Nagle, R\"odl and Schacht~\cite{Count}.
In Section~\ref{sec:extensions_count} we will deduce the extension lemma 
from the
counting lemma (which corresponds to the case when~$\Hy$ is empty).

Before this, in Section~\ref{sec:pfemb} we use the extension lemma to prove the strengthened version
of the embedding lemma mentioned before (Lemma~\ref{emblemma}). Finally, we use the embedding
lemma together with the regularity lemma for $3$-uniform 
hypergraphs due to Frankl and R\"odl~\cite{FR}
to prove Theorem~\ref{Ramseythm}. 

\section{The embedding lemma}\label{sec:stateemb}

Before we can state the embedding lemma, we first have to introduce
some notation. Given a bipartite graph $G$ with vertex classes $A$
and $B$, we denote the number of edges of~$G$ by~$e(A,B)$. The
\emph{density} of $G$ is defined to be
$$d_G(A,B):=\frac{e(A,B)}{|A||B|}.$$
We will also use $d(A,B)$ instead of $d_G(A,B)$ if this is unambiguous.
Given $0<\delta,d\le 1$, we say that $G$ is
\emph{$(d,\delta)$-regular} if for all sets $X\subseteq A$ and
$Y\subseteq B$ with $|X|\ge \delta |A|$ and $|Y|\ge \delta |B|$ we
have $(1-\delta)d<d(X,Y)<(1+\delta)d$.

Given $\delta>0$, we also say that $G$ is \emph{$\delta$-regular} if for
any $X \subseteq A$ and $Y \subseteq B$ which satisfy $|X| \geq \delta
|A|$, $|Y| \geq \delta |B|$, we have
$$|d(X,Y)-d(A,B)| \leq \delta.$$
It can easily be seen that this definition of regularity is roughly
equivalent to $(d,\delta)$-regularity.
We say that a \emph{$k$-partite graph~$P$ is $(d,\delta)$-regular}
if each of the~$\binom{k}{2}$ bipartite subgraphs forming~$P$ is
$(d,\delta)$-regular or empty.

Given a $3$-uniform hypergraph~$\G$, we denote by $|\G|$ the number
of its vertices and by~$E(\G)$ the set of its hyperedges. We write $e(\G):=|E(\G)|$.
We say that vertices $x,y\in\G$ are \emph{neighbours} if~$x$ and~$y$ lie in
a common hyperedge of~$\G$.

In order to state the embedding lemma we will now say what we mean by a `pseudo-random'
hypergraph, i.e.~we will formally define regularity of $3$-uniform hypergraphs.
Suppose we are given a 3-partite
graph~$P$ with vertex classes $V_i,V_j,V_k$ where the three bipartite
graphs forming~$P$ are denoted by $P^{ij}$, $P^{jk}$ and~$P^{ik}$.
We will often refer to such a 3-partite graph as a \emph{triad}.
We write $T(P)$ for the set of all triangles contained in $P$ and
let $t(P)$ denote the number of these triangles. Given a $3$-partite
$3$-uniform hypergraph $\G$ with the same vertex classes, we define
the \emph{density of  $P$ with respect to $\G$} by $$ d_{\G}(P):=
\begin{cases}
|E(\G)\cap T(P)|/t(P) & \text{if } t(P)>0,\\
0 & \text{otherwise.}
\end{cases}
$$
In other words, $d_{\G}(P)$ denotes the proportion of all those
triangles in~$P$ which are hyperedges of~$\G$. More generally,
suppose that we are given an $r$-tuple $\vec{Q}=(Q(1),\dots,Q(r))$
of subtriads of~$P$, where $Q(s)=Q^{ij}(s)\cup Q^{jk}(s) \cup
Q^{ik}(s)$, and $Q^{ij}(s)\subseteq P^{ij}$, $Q^{jk}(s)\subseteq
P^{jk}$, $Q^{ik}(s)\subseteq P^{ik}$ for all $s \in [r]$, where
$[r]$ denotes $\lbrace1,\dots,r\rbrace$. Put
$$ t(\vec{Q}):=\left|\bigcup_{s=1}^r T(Q(s))\right|. $$ The
\emph{density of $\vec{Q}$ with respect to $\G$} is defined to be $$
d_{\G}(\vec{Q}):= \begin{cases} |E(\G)\cap \bigcup_{s=1}^r
T(Q(s))|/t(\vec{Q}) & \text{if }
t(\vec{Q})>0,\\
0 & \text{otherwise.}
\end{cases}
$$
\noindent Note that in this definition, the sets $T(Q(s))$ of
triangles need not necessarily be disjoint. We say that a triad $P$
is \emph{$(d_3,\delta_3,r)$-regular} with respect to $\G$ if for
every $r$-tuple $\vec{Q}=(Q(1),\dots,Q(r))$ of subtriads of $P$ with
$$t(\vec{Q})\ge \delta_3\cdot t(P)$$ we have
$$|d_3-d_{\G}(\vec{Q})|<\delta_3.$$
%
We say that $P$ is \emph{$(\delta_3,r)$-regular} with respect to~$\G$
if it is $(d,\delta_3,r)$-regular for some $d$. More generally, if
$k\ge 3$, $P$ is a $k$-partite graph and~$\G$ is a $k$-partite
$3$-uniform hypergraph with the same vertex classes, we say that~$P$
is \emph{$(d_3,\delta_3,r)$-regular} with respect to~$\G$ if each of
the triads~$P'$ induced by~$P$ is either $(d_3,\delta_3,r)$-regular with
respect to~$\G$ or satisfies~$d_\G(P')=0$. 

If $\Hy$ and~$\G$ are
$k$-partite $3$-uniform hypergraphs~ with vertex classes
$X_1,\dots,X_k$ and $V_1,\dots,V_k$ respectively, and if $P$ is a
$k$-partite graph with vertex classes $V_1,\dots,V_k$, we say that
\emph{$(\G,P)$ respects the partition of $\Hy$} if, for all
$i<j<\ell$, whenever~$\Hy$ contains a hyperedge with vertices in
$X_i,X_j,X_\ell$, the hypergraph~$\G$ contains a hyperedge with
vertices in $V_i,V_j,V_\ell$ which also forms
a triangle in~$P$.

Note that if $(\G,P)$ respects the partition of $\Hy$ and $P$ is
$(d_3,\delta_3,r)$-regular with respect to~$\G$ then the triad
$P[V_i,V_j,V_\ell]$ induced by~$V_i\cup V_j\cup V_\ell$ is $(d_3,\delta_3,r)$-regular
whenever~$\Hy$ contains a hyperedge with vertices in $X_i,X_j,X_\ell$. Thus if
$P$ is also graph-regular, 
$\Hy$ has bounded maximum degree, $\delta_3\ll d_3$ and $|X_i|\le
|V_i|$ for all~$i$, then one might hope that this regularity can be
used to find an embedding of~$\Hy$ in~$\G$ (where the vertices
in~$X_i$ are represented by vertices in~$V_i$).

\begin{lemma} [Embedding lemma]\label{embcor}
Let $\Delta,k,r,n_0$ be positive integers and let
$c,d_2,d_3,\delta_2,\delta_3$ be positive constants such that
$$1/n_0\ll1/r \ll \delta_2 \ll\min\{\delta_3,d_2\}\le\delta_3\ll
d_3,1/\Delta,1/k \text{\ \ \ and\ \ \ }
c\ll d_2,d_3,1/\Delta,1/k.$$
Then the following holds for all integers $n\ge
n_0$. Suppose that $\Hy$ is a $k$-partite $3$-uniform hypergraph of
maximum degree at most~$\Delta$ with vertex classes $X_1,\dots,X_k$
such that $|X_i|\le cn$ for all $i=1,\dots,k$. Suppose that $\G$ is
a $k$-partite hypergraph with vertex classes $V_1,\dots,V_k$, which
all have size~$n$. Suppose that $P$ is a $(d_2,\delta_2)$-regular
$k$-partite graph with vertex classes $V_1,\dots,V_k$ which is
$(d_3,\delta_3,r)$-regular with respect to~$\G$, and $(\G,P)$
respects the partition of~$\Hy$. Then~$\G$ contains a copy of~$\Hy$.
\end{lemma}

Here we write $0<a_1 \ll a_2 \ll a_3$ to mean
that we can choose the constants $a_1,a_2,a_3$ from right to left. More
precisely, there are increasing functions $f$ and $g$ such that, given
$a_3$, whenever we choose some $a_2 \leq f(a_3)$ and $a_1 \leq g(a_2)$, all
calculations needed in the proof of Lemma~\ref{embcor} are valid. In order to
simplify the exposition, we will not determine these functions explicitly.
Hierarchies with more constants are defined in the obvious way.

The strategy of our proof of Lemma~\ref{embcor} is to proceed by
induction on~$|\Hy|$. So for any vertex~$h$ of~$\Hy$, let~$\Hy_h$
denote the hypergraph obtained from $\Hy$ by removing~$h$. Let $v,w$
be any vertices of~$\Hy$ forming a hyperedge with~$h$. In the
induction step, we only want to consider copies of $\Hy_h$ in $\G$
for which $vw$ is an edge of $P$ (otherwise there is clearly no
chance of using the regularity of~$\G$ to extend this copy of~$\Hy_h$
to one of~$\Hy$). This
motivates the following definition. A \emph{complex~$\Hy$} consists
of vertices, edges and hyperedges such that the set of edges is a
subset of the set of unordered pairs of vertices and the set of
hyperedges is a subset of the set of unordered triples of vertices.
Moreover, each pair of vertices of~$\Hy$ lying in a common hyperedge
has to form an edge of~$\Hy$. Thus we can make every $3$-uniform hypergraph~$\Hy$
into a complex by adding an edge between every pair of vertices that lies in
a common hyperedge of~$\Hy$. We will often denote this complex by~$\Hy$ again.

Instead of Lemma~\ref{embcor}, we will prove an embedding lemma for complexes.
In order to state it, we need to introduce some more notation.
Given a complex $\Hy$, we let $V(\Hy)$
denote the set of its vertices, we write $E_2(\Hy)$ for the set
of its edges and $E_3(\Hy)$ for the set of its hyperedges. Note
that each hyperedge of a complex~$\Hy$ forms a triangle in the
underlying graph (whose vertex set is~$V(\Hy)$ and whose set of
edges is~$E_2(\Hy)$). We set $|\Hy|:=|V(\Hy)|$, and
$e_i(\Hy):=|E_i(\Hy)|$ for $i=2,3$. We say that a complex~$\Hy$ is
\emph{$k$-partite} if its underlying graph is $k$-partite. The \emph{degree of
a vertex~$x$ in a complex~$\Hy$} is the maximum%
     \COMMENT{That's a new definition which simplifies some of the calculations
later on.}
of the degree of~$x$ in
the underlying graph and its degree in the underlying hypergraph (whose
vertex set is~$V(\Hy)$ and whose set of hyperedges is~$E_3(\Hy)$).
The \emph{maximum degree of~$\Hy$} is then defined in the obvious way. 
We say that vertices~$x$ and $y$ are \emph{neighbours in~$\Hy$} if they
are neighbours in the underlying graph. Subcomplexes
of~$\Hy$ and subcomplexes induced by some vertex set~$X\subseteq
V(\Hy)$ are defined in the natural way. Also, the symbol $K_k^{(3)}$
will denote either the complete complex or the complete $3$-uniform
hypergraph on $k$ vertices. It will be clear from the context which
of the two is intended. Note that the complete complex~$K_1^{(3)}$ is
just a vertex and~$K_2^{(3)}$ consists of two vertices joined by an edge.

Given $k$-partite complexes~$\Hy$ and~$\G$ with vertex classes
$X_1,\dots,X_k$ and $V_1,\dots,V_k$, we say that \emph{$\G$ respects
the partition of $\Hy$} if it statisfies the following two
properties. Firstly, for all $i<j<\ell$, the complex~$\G$ contains
a hyperedge with vertices in $V_i,V_j,V_\ell$ whenever~$\Hy$
contains a hyperedge with vertices in $X_i,X_j,X_\ell$. Secondly,
for all $i<j$, the complex~$\G$ contains an edge between~$V_i$
and~$V_j$ whenever~$\Hy$ contains an edge between~$X_i$ and~$X_j$.

We say that a $k$-partite complex $\G$ is
\emph{$(d_3,\delta_3,d_2,\delta_2,r)$-regular} if its underlying
graph~$P$ is $(d_2,\delta_2)$-regular and~$P$ is
$(d_3,\delta_3,r)$-regular with respect to the underlying
hypergraph of~$\G$.

Suppose that we have $k$-partite complexes~$\Hy$
and~$\G$ with vertex classes $X_1,\dots,X_k$ and $V_1,\dots,V_k$
respectively. A \emph{labelled partition-respecting copy of $\Hy$ in
$\G$} is a labelled subcomplex of~$\G$ which is isomorphic to $\Hy$
such that the corresponding isomorphism maps~$X_i$ to a subset
of~$V_i$. This definition naturally extends to labelled
partition-respecting copies of subcomplexes~$\Hy'$ of~$\Hy$ in~$\G$.
Given any subcomplex $\Hy'$ of $\Hy$, we write $|\Hy'|_\G$ for the
number of labelled partition-respecting copies of $\Hy'$ in~$\G$.

Instead of Lemma~\ref{embcor} we will prove the following result
which implies it immediately.

\begin{lemma}[Embedding lemma for complexes]\label{emblemma}
Let $\Delta,k,r,n_0$ be positive
integers and let $c,\alpha,d_2,d_3,\delta_2,\delta_3$ be positive
constants such that
$$1/n_0\ll 1/r \ll \delta_2 \ll\min\{\delta_3,d_2\}\le\delta_3\ll
\alpha \ll d_3,1/\Delta,1/k \text{\ \ \ and\ \ \ }
c\ll \alpha, d_2.$$
Then the following holds for all integers
$n\ge n_0$. Suppose that $\Hy$ is a $k$-partite complex of maximum
degree at most~$\Delta$ with vertex classes $X_1,\dots,X_k$ such
that $|X_i|\le cn$ for all $i=1,\dots,k$. Suppose also that $\G$ is
a $k$-partite $(d_3,\delta_3,d_2,\delta_2,r)$-regular complex with
vertex classes $V_1,\dots,V_k$, all of size $n$, which respects the
partition of~$\Hy$. Then for every vertex $h\in \Hy$ we have that
$$
|\Hy|_{\G} \geq (1 - \alpha) n d_2^{e_2(\Hy)-e_2(\Hy_h)}
d_3^{e_3(\Hy)-e_3(\Hy_h)} |\Hy_h|_{\G},
$$
where $\Hy_h$ denotes the induced subcomplex of~$\Hy$ obtained by
removing $h$. In particular, $\G$ contains at least
$((1-\alpha)n)^{|\Hy|} d_2^{e_2(\Hy)}d_3^{e_3(\Hy)}$ labelled
partition-respecting copies of $\Hy$.
\end{lemma}

Note that we would expect almost
$n d_2^{e_2(\Hy)-e_2(\Hy_h)}d_3^{e_3(\Hy)-e_3(\Hy_h)} |\Hy_h|_{\G}$
labelled partition-respecting copies of~$\Hy$ if~$\G$ were a random complex.
As indicated above, we will prove
Lemma~\ref{emblemma} by induction on~$|\Hy|$. In the induction step, it will be extremely useful
to assume the existence of the expected number of copies of any proper
subcomplex~$\Hy'$ of~$\Hy$
in~$\G$ and not just the existence of one such copy.

\section{Tools}\label{sec:tools}

In our proof of Lemma~\ref{emblemma} we will use the so-called
counting lemma.

\begin{lemma}[Counting lemma]\label{countinglemma}
Let $k,r,t,n_0$ be positive integers and let
$\beta,d_2,d_3,\delta_2,\delta_3$ be positive constants such that
$$1/n_0\ll 1/r \ll \delta_2 \ll\min\{\delta_3,d_2\}\le\delta_3\ll
\beta,d_3,1/k,1/t.$$ Then the following holds for all integers
$n\ge n_0$. Suppose that $\Hy$ is a $k$-partite complex on $t$
vertices with vertex classes $X_1,\dots,X_k$. Suppose also that $\G$
is a $k$-partite $(d_3,\delta_3,d_2,\delta_2,r)$-regular complex
with vertex classes $V_1,\dots,V_k$, all of size $n$, which respects
the partition of~$\Hy$. Then~$\G$ contains $$(1\pm\beta)
n^{t} d_2^{e_2(\Hy)}d_3^{e_3(\Hy)} $$ labelled partition-respecting
copies of~$\Hy$.
\end{lemma}

The lower bound in Lemma~\ref{countinglemma} for~$K_k^{(3)}$'s was proved by
Nagle and~R\"odl~\cite{NR03}
(a short proof was given later in~\cite{Nrs}). 
Nagle, R\"odl and Schacht~\cite{Count} generalized this
lower bound to arbitrary $k$-uniform complexes (Lemma~\ref{counting2} in
Section~\ref{uppercount}).
In a slightly different setup, this
was also proved independently by Gowers~\cite{Gowers}.%
     \COMMENT{who has an upper bound as well}
The upper bound in
Lemma~\ref{countinglemma} can easily be derived from the lower bound. This was done
for~$K_k^{(3)}$'s in~\cite{NR03}.
In Section~\ref{uppercount} we show how one can derive
Lemma~\ref{countinglemma} from Lemma~\ref{counting2}.

Note that Lemma~\ref{emblemma} is a generalization of the lower bound in
Lemma~\ref{countinglemma}.
As a special case, Lemma~\ref{countinglemma} includes the counting lemma
for graphs, which is an easy consequence of the definition of graph regularity.

The following result is another strengthening of
Lemma~\ref{countinglemma}. We will need it in the proof of Lemma~\ref{emblemma}.
Given complexes $\Hy\subseteq \Hy'$ such
that $\Hy$ is induced, it states that $\G$ not only contains about
the expected number of copies of~$\Hy'$, but also that almost all
copies of~$\Hy$ in~$\G$ are extendible to about the expected number of
copies of~$\Hy'$. The special case when~$\Hy$ is a hyperedge was proved earlier by
Haxell, Nagle and R\"odl~\cite{haxell3}.

\begin{lemma}[Extension lemma]\label{extensions_count}
Let $k,r,t,t',n_0$ be positive integers, where $t<t'$, and let
$\beta,d_2,d_3,\delta_2,\delta_3$ be positive constants such that
$$1/n_0\ll 1/r \ll \delta_2 \ll\min\{\delta_3,d_2\}\le\delta_3\ll
\beta,d_3,1/k,1/t'.$$ Then the following holds for all integers $n\ge n_0$.
Suppose that $\Hy'$ is
a $k$-partite complex on $t'$ vertices with vertex classes
$X_1,\dots,X_k$ and let~$\Hy$ be an induced subcomplex of~$\Hy'$ on~$t$
vertices. Suppose also that $\G$ is a $k$-partite
$(d_3,\delta_3,d_2,\delta_2,r)$-regular complex with vertex classes
$V_1,\dots,V_k$, all of size $n$, which respects the partition of
$\Hy'$. Then all but at most $\beta |\Hy|_{\G}$  labelled
partition-respecting copies of~$\Hy$ in~$\G$ are extendible into 
$$(1 \pm \beta) n^{t'-t} d_2^{e_2(\Hy')-e_2(\Hy)}d_3^{e_3(\Hy')-e_3(\Hy)}
$$
labelled partition-respecting copies of~$\Hy'$ in~$\G$.
\end{lemma}

Lemmas~\ref{countinglemma} and~\ref{extensions_count} differ from
Lemma~\ref{emblemma} in that the positions of~$t$
and~$t'$ in the hierarchy mean we can only look at complexes~$\Hy,\Hy'$
of bounded size. In particular, in the proof of Lemma~\ref{emblemma}
we will apply these lemmas to complexes whose order is some function of~$\Delta$
and so does not depend on~$n$. Lemma~\ref{extensions_count} will
be proved in Section~\ref{sec:extensions_count}.


\section{Proof of the embedding lemma}\label{sec:pfemb}

Throughout this section, if we refer to
a copy of a certain subcomplex~$\Hy'$ of~$\Hy$ in~$\G$ we mean that
this copy is labelled and partition-respecting without mentioning
this explicitly. We usually denote such a copy of~$\Hy'$ by~$H'$
(i.e.~by the corresponding roman letter).
We fix new constants $\beta$ and
$\delta'_2$ such that
$$
\delta_2\ll\delta'_2\ll d_2,d_3,1/\Delta$$ 
and 
$$c,\delta'_2,\delta_3\ll \beta\ll \alpha.$$ 
As mentioned earlier, we will prove Lemma~\ref{emblemma} by
induction
on~$|\Hy|$. 
We first show that we may assume that the component $\C_1$ of $\Hy$
containing~$h$ satisfies $|\C_1| > \Delta^5$.
So suppose this is not the case and let $\C_2:=\Hy-\C_1$.
Every copy of $\Hy$ can be obtained by first choosing a copy $C_2$
of $\C_2$ and then choosing a copy $C_1$ of $\C_1$ which is disjoint
from $C_2$. Thus
\begin{equation} \label{eqlowerC}
|\Hy|_\G = \sum_{C_2 \in \G} |\C_1|_{\G-C_2}
\ge \sum_{C_2 \in \G} \frac{(1-c)^{\Delta^5}(1-\beta)}{1+\beta} |\C_1|_\G
\ge (1-3\beta)|\C_1|_\G |\C_2|_\G.
\end{equation}
Here we applied the counting lemma (Lemma~\ref{countinglemma}) 
in $\G-\C_2$ and in $\G$ to obtain the first inequality.
On the other hand, for $\Hy_h:=\Hy-h$ we have
\begin{equation} \label{equpperC}
|\Hy_h|_\G \le | \C_1-h|_\G |\C_2|_\G
\le \frac{(1+\beta)
|\C_1|_\G |\C_2|_\G}{(1-\beta) d_2^{e_2(\C_1)-e_2(\C_1-h)}d_3^{e_3(\C_1)-e_3(\C_1-h)}n},
\end{equation}
where the second inequality follows from the application of the counting lemma to $\C_1$ and $\C_1 - h$.
Combining~(\ref{eqlowerC}) and~(\ref{equpperC}) gives the result claimed above.
Note that in particular, this deals with the start of the induction.
So we may assume that $|\Hy|> \Delta^5$ and that  Lemma~\ref{emblemma} 
holds for all complexes with
fewer than $|\Hy|$ vertices.
Also, the above assumption on~$\C_1$ together with the fact that~$\Hy$ has maximum degree 
$\Delta$ implies that the set of all those vertices of $\Hy$ which (in the underlying graph)
have distance exactly $3$ to~$h$ is nonempty. This will be convenient later on.

For induced subcomplexes $\Hy''\subseteq \Hy'\subseteq\Hy$ and a
copy~$H''$ of~$\Hy''$ in~$\G$, we denote by $|H''\to \Hy'|_\G$ the
number of copies of~$\Hy'$ in~$\G$ which extend~$H''$. We set
$$\overline{|\Hy''\to
\Hy'|}:=d_2^{e_2(\Hy')-e_2(\Hy'')}d_3^{e_3(\Hy')-e_3(\Hy'') }
n^{|\Hy'|-|\Hy''|}.$$ Thus $\overline{|\Hy''\to \Hy'|}$ is roughly
the expected number of ways a copy of~$\Hy''$ in~$\G$ could be
extended to a copy of~$\Hy'$ if~$\G$ were a random complex.

Given the vertex $h\in\Hy$ as in Lemma~\ref{emblemma}, we write~$\cn$
for the subcomplex of~$\Hy$ induced by all the
neighbours of~$h$ in~$\Hy$. We write $\B$ for the
subcomplex of~$\Hy$ induced by~$V(\cn)\cup\{h\}$.
We call a copy~$N_h$ of~$\cn$ in $\G$ \emph{typical}
if~$N_h$ can be extended to at least $(1- \beta)\overline{|\cn
\rightarrow \B|}$ copies of $\B$.
If we knew that every copy of~$\cn$ in~$\G$ were typical, then the induction
step would follow immediately since this would imply that~$|\Hy|_\G$ is
roughly %
     \COMMENT{We don't have equality instead of just `roughly' since
we might choose vs more than once.}
$$
\sum_{N_h\in \G} |N_h\to\Hy_h|_\G |N_h\to \B|_\G 
\ge (1-\beta)\overline{|\cn \rightarrow \B|} \sum_{N_h\in \G} |N_h\to\Hy_h|_\G
= (1-\beta)\overline{|\cn \rightarrow \B|}|\Hy_h|_\G.
$$
Indeed, this would hold since each copy of $\Hy$ in~$\G$ can be obtained by first
choosing a copy~$N_h$ of $\cn$, then extending~$N_h$ to some copy of $\Hy_h$
and then extending~$N_h$ to a copy of~$\B$.
However, the extension lemma (Lemma~\ref{extensions_count}) only implies that
almost all copies of~$\cn$ are typical, which makes things more complicated.
So let~$\Typ$ denote the set of
all typical copies of~$\cn$ in~$\G$ and $\Atyp$ the set of all other copies.
Lemma~\ref{extensions_count} implies that
\begin{equation}\label{eqtyp} |\Typ|\geq (1-\beta)|\cn |_{\G}.
\end{equation}

We now define an analogous set where we refer to the underlying
graph~$P$ of~$\G$ instead of~$\G$ itself. More precisely, we call a
copy~$N_h$ of~$\cn$ \emph{useful} if the following holds: Let
$x_1,\dots,x_\ell$ be any distinct vertices of $\n$ and
let $x_1',\dots,x_\ell'$ be the corresponding vertices in~$\cn$.
If a vertex class~$X_i$ contains a common neighbour of $x_1',\dots,x_\ell'$,
then in the underlying graph~$P$ the common neighbourhood
of $x_1,\dots , x_\ell$ in $V_i$ has size
$(1\pm \delta_2)^\ell d_2^\ell n$. We denote by~$\Usef$ the set of all these
copies of~$\cn$.

We will now show that almost all copies of~$\cn$ in~$\G$ are useful.
First recall that since~$\G$ respects the partition of~$\Hy$, the bipartite graphs
$P[V_i,V_j]$ are $(d_2,\delta_2)$ regular whenever~$\Hy$ contains an edge between
$X_i$ and~$X_j$. Together with the fact that $|\cn| \le \Delta$ this shows
at most at $2\Delta^2 2^\Delta\delta_2 n^{|\cn|}$ of the $|\cn|$-tuples of
vertices in~$\G$ do not satisfy the above neighbourhood condition in some of the relevant
vertex classes~$V_i$. Indeed, to see this first note that the graph regularity
implies that each vertex class contains at most $2\delta_2n$
vertices having degree $\neq (1\pm \delta_2)d_2|A|$ in any given sufficiently large subset~$A$
of~$V_i$. Thus the number of
$\ell$-tuples $x_1,\dots,x_\ell$ of vertices in~$\G$ whose common neighbourhood in~$V_i$
has size $\neq (1\pm \delta_2)^\ell d_2^\ell n$ is at most~$2\ell\delta_2 n^\ell$.
Given $x_1',\dots,x_\ell'$, there are at most~$\Delta$ choices for~$V_i$.
The bound now follows since there are at most $2^\Delta$ choices for $\{x_1',\dots,x_\ell'\}$.

On the other hand, Lemma~\ref{countinglemma} implies that
$|\cn|_\G \ge \frac{1}{2}(d_2d_3)^{\Delta^2}n^{|\cn|}$. Altogether this shows that
\begin{equation}\label{equsef}
|\Usef|\ge |\cn |_{\G} -2\Delta^2 2^\Delta \delta_2 n^{|\cn|}
\geq (1-\delta'_2)|\cn |_{\G}.
\end{equation}

Recall that each copy of $\Hy$ in~$\G$ can be obtained by first
choosing a copy~$N_h$ of $\cn$, then extending~$N_h$ to some
copy~$H_h$ of $\Hy_h$ and then extending~$N_h$ to a copy of~$\B$. In
the final step we have to choose a vertex $x\in \G$ which can play
the role of~$h$. If~$N_h$ is typical then there are at least $(1-
\beta)\overline{|\cn \rightarrow \B|}$ possible choices for~$x$.
However, we have to make sure that~$x$ does not already lie
in~$H_h$. The latter condition excludes at most $cn\le
\beta\overline{|\cn \rightarrow \B|}$ of the possible choices
for~$x$. So altogether we have that
\begin{eqnarray}\label{first}
|\Hy|_{\G} &\geq & (1-2\beta)\overline{|\cn\to \B|}\sum_{N_h\in
\Typ}
|N_h\to\Hy_h |_\G\nonumber \\
&= & (1-2\beta )\overline{|\cn \to \B|} \left( \sum_{N_h\in \G}
|N_h\to\Hy_h|_\G
- \sum_{N_h\in \Atyp} |N_h\to\Hy_h|_\G\right)\nonumber \\
&\ge & (1-2\beta )\overline{|\cn \to \B|} \left( |\Hy_h|_\G -
\sum_{N_h\in \Atyp\cap \Usef} |N_h\to\Hy_h|_\G -
\sum_{N_h\notin\Usef} |N_h\to\Hy_h|_\G\right).
\end{eqnarray}

So our aim now is to prove that each of the last two sums
in~(\ref{first}) contributes no more than a small proportion of
$|\Hy_h|_{\G}$. More precisely, we will show that
\begin{equation} \label{two_sums}
\sum_{N_h\in \Atyp\cap \Usef} |N_h\to\Hy_h|_\G +
\sum_{N_h\notin\Usef} |N_h\to\Hy_h|_\G\le \beta^{1/2} |\Hy_h|_{\G}.
\end{equation}
Since $\beta\ll \alpha$ this then proves the induction step. To
prove~(\ref{two_sums}), we bound both sums separately. 
In both cases, we bound $|N_h\to\Hy_h|_\G$ in terms of its average
value  
$$
\frac{1}{|\cn|_\G} \sum_{N_h \in \G} |N_h\to\Hy_h|_\G   = \frac{|\Hy_h|_\G }{|\cn|_\G}.
$$ 
Our upper
bound for the first sum in~(\ref{two_sums}) will follow easily from the next claim.

\medskip

\noindent \textbf{Claim 1.} \emph{Every useful copy~$N_h$
of~$\cn$ in~$\G$ satisfies \[|\n \to \Hy_h|_{\G} \leq
 \frac{12}{d_3^{2\Delta^3}}
 \frac{|\Hy_h|_{\G}}{|\cn|_{\G}}. \]}
\smallskip

\noindent To prove this claim, we fix any useful copy~$\n$ of~$\cn$.
We let $\cn^*$ be the subcomplex of~$\Hy$ induced by the vertices
which have distance 2 to the vertex set of $\cn$ in the underlying
graph. Recall that our assumption at the beginning of the proof of
the lemma implies that $\cn^*$ is nonempty.
Moreover, $|\cn^*| \le \Delta^2 |\cn| \le \Delta^3$.
Let $\F' \subseteq \Hy$ be the subcomplex of~$\Hy$ that is induced
by $V(\cn)\cup V(\cn^*)$ and all the vertices in the first
neighbourhood of~$\cn$ in~$\Hy_h$ (see Figure~1).
\begin{figure} \footnotesize
\centering
\psfrag{1}[][]{$h$}
\psfrag{2}[][]{$\ \cn$}
\psfrag{3}[][]{$\ \cn^*$}
\psfrag{4}[][]{$\Hy_h^{*}$}
\psfrag{5}[][]{$\F'$}
\psfrag{6}[][]{$\Hy_h$}
\includegraphics[scale=0.5]{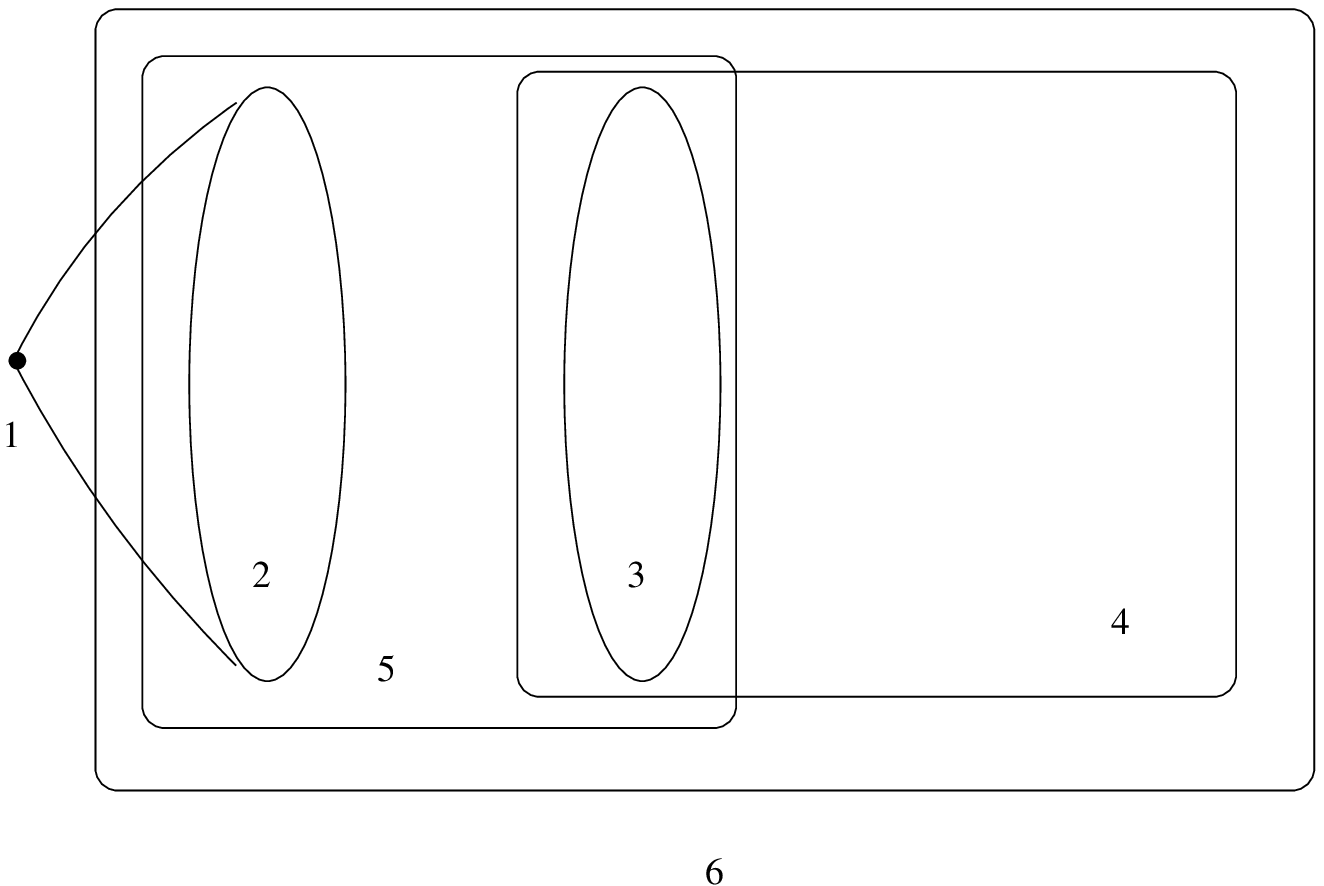}
\caption{The hypergraph $\Hy$}
\end{figure}%
So $h \not \in V(\F')$. Let $\F$ denote the underlying graph of~$\F'$. Given a
copy $\n^*$ of~$\cn^*$ in~$\G$, we denote by $|\n, \n^*
\stackrel{P}{\to} \F|_{\G}$ the number of ways the underlying graphs
of~$\n$ and~$\n^*$ can be extended into a copy of $\F$ (within the
graph~$P$). Similarly, we set
\begin{equation}\label{eqoverF}
\overline{|\cn, \cn^*\stackrel{P}{\to} \F|}:= d_2^{e_2(\F')-e_2(\cn)-e_2(\cn^*)}
n^{|\F'|-|\cn|-|\cn^*|}.
\end{equation}
Thus $\overline{|\cn, \cn^*\stackrel{P}{\to} \F|}$ is roughly the expected
number of ways the underlying graphs of disjoint copies of~$\cn$ and~$\cn^*$
can be extended into a copy of the graph $\F$, if~$\G$ were a random complex.

We define a copy~$\n^*$ of $\cn^*$ 
to be \emph{useful with respect to~$\n$} if it is disjoint from $\n$ and
if the following holds. Let $x_1,\dots,x_\ell$ and $y_1,\dots,y_{\ell^*}$ be any distinct
vertices of~$\n$ and~$\n^*$ respectively. Let $x'_1,\dots,x'_\ell$ and $y'_1,\dots,y'_{\ell^*}$
denote the corresponding vertices in~$\cn$ and~$\cn^*$. If a vertex
class~$X_i$ of~$\Hy$ contains
a common neighbour of $x'_1,\dots,x'_\ell,y'_1,\dots,y'_{\ell^*}$ in~$\F-V(\cn\cup\cn^*)$
then in the underlying graph~$P$ the common neighbourhood of
$x_1,\dots,x_\ell,y_1,\dots,y_{\ell^*}$
in~$V_i$ has size $(1\pm \delta_2)^{\ell+\ell^*}d_2^{\ell+\ell^*}n$.
We denote the set of all such copies of~$\cn^*$ in~$\G$ by~$\Usef^*(\n)$.
Using the fact that~$\n$ is useful, similarly as in~(\ref{equsef}) one can show that
\begin{equation}\label{suitable_copies_1}
|\Usef^*(\n)|\geq (1-\delta_2')|\cn^* |_{\G}.
\end{equation}
(Note that the condition that a useful copy of $\cn^*$ has to be disjoint from 
$\n$ does not affect the calculation significantly.)
Moreover, since all the bipartite subgraphs forming~$P$ are
$(d_2,\delta_2)$-regular or empty, we see that every~$\n^* \in \Usef^*(\n)$ satisfies
\begin{equation}\label{suitextension} |\n, \n^* \stackrel{P}{\to}
\F|_{\G} \leq 2 \overline{|\cn, \cn^* \stackrel{P}{\to} \F|}.
\end{equation}
Indeed, let $\F^*:=\F-V(\cn\cup\cn^*)$ and let $w_1,\dots,w_p$ denote the vertices
of~$\F^*$. 
Let $N'(w_i)$ be the neighbourhood of~$w_i$ in~$V(\cn\cup\cn^*)$ in the
graph~$\F$. Let~$W_i$
denote the set of candidates for~$w_i$ inside the vertex class of~$\G$ which we aim to
embed~$w_i$ into. Thus~$W_i$ consists of all those vertices in that class which are joined
to all the vertices in $\n\cup\n^*$ corresponding to~$N'(w_i)$.
The usefulness of~$\n$ and~$\n^*$ implies that $|W_i|=((1\pm \delta_2)d_2)^{|N'(w_i)|}n$.
In particular, the subgraph of~$P$ induced by the~$W_i$'s is still regular.
So the counting lemma for graphs implies that the number of copies of~$\F^*$
induced by the~$W_i$'s is at most
$$
\frac{3}{2} d_2^{e(\F^*)}\prod_{i=1}^p |W_i|
\stackrel{(\ref{eqoverF})}{\le} 2 \overline{|\cn, \cn^* \stackrel{P}{\to} \F|},
$$
as required.

Let $\Hy_h^*$ denote the subcomplex of~$\Hy$ obtained by
deleting~$h$ as well as all the vertices in $\F'-\cn^*$. Then any
copy of~$\Hy_h$ extending~$\n$ can be obtained by first choosing a
copy~$\n^*$ of~$\cn^*$, then extending this copy to a copy~$H^*_h$
of~$\Hy^*_h$, and then extending the pair~$\n$, $\n^*$ into a copy
of~$\F'$ (which avoids~$H^*_h$). Clearly, there are at most $|\n,
\n^* \stackrel{P}{\to} \F|_{\G}$ ways to choose an extension
of~$\n$, $\n^*$ into a copy of~$\F'$. 
(Using the latter bound means that we are disregarding any hyperedges
of $\Hy_h$ in $E_3(\F') \setminus E_3(\cn \cup \cn^*)$.
This is the reason for the error term involving $d_3$ in the statement
of Claim~1.) Thus
\begin{eqnarray}
\label{first_upper_bound} \lefteqn{|\n \to \Hy_h|_{\G} \leq} \nonumber \\
& &\sum_{\n^* \in \Usef^*(\n)} |\n, \n^* \stackrel{P}{\to} \F|_{\G}
|\n^* \to \Hy_h^*|_{\G} + \sum_{\n^* \not \in \Usef^*(\n)} |\n, \n^*
\stackrel{P}{\to} \F|_{\G} |\n^* \to \Hy_h^*|_{\G}.
\end{eqnarray}
The first sum in~(\ref{first_upper_bound}) can be bounded by
\begin{eqnarray}\label{eqsuit}
\sum_{\n^* \in \Usef^*(\n)} |\n, \n^* \stackrel{P}{\to} \F|_{\G}
|\n^* \to \Hy_h^*|_{\G} &\stackrel{(\ref{suitextension})}{\le }& 2
\overline{|\cn, \cn^* \stackrel{P}{\to} \F|} \sum_{\n^* \in
\Usef^*(\n)} |\n^* \to \Hy_h^*|_{\G}\nonumber \\ &\leq & 2
\overline{|\cn, \cn^* \stackrel{P}{\to} F|}  |\Hy_h^*|_{\G}.
\end{eqnarray}
To bound the second sum in (\ref{first_upper_bound}),
let~$\Hy'_h:=\Hy^*_h-\cn^*$. Then clearly
\begin{equation}\label{eqHyh} |\n^* \to \Hy_h^*|_{\G}\le
|\Hy'_h|_\G.
\end{equation}
We shall estimate $|\Hy'_h|_{\G}$ in relation to  $|\Hy_h^*|_{\G}$.
Let $s:=|\cn^*|=|\Hy_h^*|-|\Hy'_h|$ and suppose that $w_1,\ldots,
w_s$ are the vertices in~$\Hy_h^*-\Hy'_h=\cn^*$. For all $i=1,\ldots, s$,
we let $\N'_i$ be the subcomplex induced by the neighbourhood of $w_i$ in
$\Hy_h^*-\{w_{i+1},\dots, w_s\}$. Let $\B'_i$ be the subcomplex
of~$\Hy$ induced by~$w_i$ and
the vertices in~$\N'_i$. Then our induction hypothesis implies that%
     \COMMENT{We don't write just $s$ for $|\Hy_h^*|-|\Hy'_h|$ in the
inequality below since we will need to refer to it later on, and it
will be clearer in this form.}
\begin{eqnarray}\label{smalltolargefirst} |\Hy_h^* |_{\G} & \geq &
(1-\alpha)^{|\Hy_h^*|-|\Hy'_h|}
\left(\prod_{i=1}^{|\Hy_h^*|-|\Hy'_h|}
\overline{|\N'_i \to \B'_i|}\right)|\Hy'_h|_{\G}\nonumber \\
& = &((1-\alpha)n)^{|\Hy_h^*|-|\Hy'_h|}
d_2^{e_2(\Hy_h^*)-e_2(\Hy'_h)}d_3^{e_3(\Hy_h^*)-e_3(\Hy'_h)}|\Hy'_h|_{\G}.
\end{eqnarray}
(In fact, one reason for our choice of the induction hypothesis is that it
allows us to relate $|\Hy_h^*|_\G$ to $|\Hy_h'|_\G$ as in~(\ref{smalltolargefirst}).)
But our assumption on the maximum degree of $\Hy$ implies that
$e_i(\Hy_h^*)-e_i(\Hy'_h) \le \Delta |\cn^*| \le \Delta^4$ for $i=2,3$. So
$$
|\Hy_h^* |_{\G}  \ge \frac{1}{2} n^{|\Hy_h^*|-|\Hy'_h|}
(d_2d_3)^{\Delta^4}|\Hy'_h|_{\G}
$$
as $\alpha\ll 1/\Delta$.
Since $|\Hy_h^*|-|\Hy'_h|=|\cn^*|$ the last inequality
together with~(\ref{suitable_copies_1}), (\ref{eqHyh}) and the fact
that $\delta'_2\ll d_2,d_3,1/\Delta$ imply that
\begin{equation}\label{second_term_second}
\sum_{\n^* \notin \Usef^*(\n)} |\n^* \to \Hy_h^*|_{\G} 
\le \delta_2' |\cn^*|_\G |\Hy_h'|_\G 
\le |\Hy_h^*|_{\G}
\frac{2\delta_2'|\cn^*|_\G}{(d_2d_3)^{\Delta^4}n^{|\cn^*|}}
\leq \sqrt{\delta'_2} |\Hy_h^*|_{\G}.
\end{equation}
In the final inequality, we also used the crude bound $|\cn^*|_\G \le
n^{|\cn^*|}$. We can now bound the second sum
in~(\ref{first_upper_bound}) by
\begin{eqnarray}\label{second_term_first}
\sum_{\n^*\notin \Usef^*(\n)} |\n, \n^* \stackrel{P}{\to} \F|_{\G}
|\n^* \to \Hy_h^*|_{\G} &\leq & n^{|\F| - |\cn|-|\cn^*|} \sum_{\n^*
\notin \Usef^*(\n)} |\n^* \to \Hy_h^*|_{\G}\nonumber\\ &
\stackrel{(\ref{second_term_second})}{\le} & \overline{|\cn, \cn^*
\stackrel{P}{\to} \F|} |\Hy_h^*|_{\G} \frac{\sqrt{\delta'_2} n^{|\F|
-
|\cn| - |\cn^*|}}{\overline{|\cn, \cn^* \stackrel{P}{\to} \F|}}\nonumber\\
& \stackrel{(\ref{eqoverF})}{\le} &  \overline{|\cn, \cn^*
\stackrel{P}{\to} \F|} |\Hy_h^*|_{\G}.
\end{eqnarray}
Indeed, to see
the last inequality recall that $\delta'_2\ll d_2,1/\Delta$.
Inequalities~(\ref{first_upper_bound}), (\ref{eqsuit})
and~(\ref{second_term_first}) together now imply that
\begin{equation}\label{eqntoHyh}
|\n \to \Hy_h|_{\G} \leq 3 \overline{|\cn, \cn^* \stackrel{P}{\to} \F|}
|\Hy_h^*|_{\G}.
\end{equation}
Similarly as in~(\ref{smalltolargefirst}) one can use the induction
hypothesis to show that
\begin{eqnarray*} |\Hy_h |_{\G} & \ge &
((1-\alpha)n)^{|\Hy_h|-|\Hy_h^*|}
d_2^{e_2(\Hy_h)-e_2(\Hy_h^*)}d_3^{e_3(\Hy_h)-e_3(\Hy_h^*)}|\Hy_h^*|_{\G}\\
& \ge  & \frac{1}{2} n^{|\F'|-|\cn^*|}
d_2^{e_2(\Hy_h)-e_2(\Hy_h^*)}d_3^{e_3(\Hy_h)-e_3(\Hy_h^*)}|\Hy_h^*|_{\G}.
\end{eqnarray*}
Observe that $e_i(\Hy_h)-e_i(\Hy_h^*) = e_i(\F')-e_i(\cn^*)$ for
$i=2,3$. Moreover, the counting lemma for graphs implies that
$|\cn|_\G \le 2 n^{|\cn|}d_2^{e_2(\cn)}$.
Together with~(\ref{eqoverF}) and (\ref{eqntoHyh}) this shows
that
\begin{eqnarray*}
|\n \to \Hy_h|_{\G} & \leq &
\frac{6d_2^{e_2(\F')-e_2(\cn)-e_2(\cn^*)}n^{|\F'|-|\cn|-|\cn^*|}|\cn|_\G}{n^{|\F'|-|\cn^*|}
d_2^{e_2(\F')-e_2(\cn^*)} d_3^{e_3(\F')-e_3(\cn^*)}}\cdot
\frac{|\Hy_h|_\G}{|\cn|_\G} \le
\frac{12|\Hy_h|_\G}{d_3^{e_3(\F')-e_3(\cn^*)}|\cn|_\G}.
\end{eqnarray*}
But $e_3(\F')-e_3(\cn^*)\le \Delta |\F - \cn^*| \le 2\Delta^3$.
This completes the proof of Claim~1.

\medskip

\noindent In order to give an upper bound on the second sum
in~(\ref{two_sums}) we will need the following claim.

\medskip

\noindent \textbf{Claim 2.} \emph{Every copy~$N_h$ of~$\cn$ satisfies
\[|N_h \to \Hy_h|_{\G} \leq 
\frac{2}{(d_2d_3)^{\Delta^{2}}}
\frac{|\Hy_h|_{\G}}{|\cn|_{\G}}. \] }

\smallskip

\noindent Let $\Hy_h^-:=\Hy_h-\cn$. Then, very crudely,
\begin{equation}\label{crudeNhtoHy}
|N_h \to \Hy_h|_{\G} \le |\Hy_h^-|_\G.
\end{equation}
But similarly as in~(\ref{smalltolargefirst})
we have that
\begin{eqnarray*} |\Hy_h
|_{\G} & \ge & ((1-\alpha)n)^{|\cn|}
d_2^{e_2(\Hy_h)-e_2(\Hy_h^-)}d_3^{e_3(\Hy_h)-e_3(\Hy_h^-)}|\Hy_h^-|_{\G}\\
& \ge & \frac{1}{2}|\cn|_\G
d_2^{e_2(\Hy_h)-e_2(\Hy_h^-)}d_3^{e_3(\Hy_h)-e_3(\Hy_h^-)}|\Hy_h^-|_{\G}\\
& \ge & \frac{1}{2}|\cn|_\G
d_2^{\Delta^2}d_3^{\Delta^2}|\Hy_h^-|_{\G}.
\end{eqnarray*}
In the final line we used the fact that $|\Hy_h-\Hy_h^-|=|\cn| \le
\Delta$. Together with~(\ref{crudeNhtoHy}) this implies Claim~2.

\medskip

\noindent Claims~1 and~2 now immediately imply~(\ref{two_sums}).
Indeed, using~(\ref{eqtyp}) and~(\ref{equsef}) and the facts that
$\delta'_2\ll d_2,d_3,1/\Delta$ and $\delta'_2\ll \beta\ll
d_3,1/\Delta$ we see that $$\sum_{N_h\in \Atyp\cap \Usef}
|N_h\to\Hy_h|_\G + \sum_{N_h\notin\Usef} |N_h\to\Hy_h|_\G\le
\left(\frac{12
\beta}{d_3^{2\Delta^3}}+\frac{2\delta'_2}{(d_2d_3)^{\Delta^2}}
\right)|\Hy_h|_{\G} \le \beta^{1/2} |\Hy_h|_{\G},$$ as required.
This completes the proof of Lemma~\ref{emblemma}.



\section{Proof of Lemma~\ref{countinglemma}} \label{uppercount}

In this section, we indicate how Lemma~\ref{countinglemma}
follows easily from the version of the counting lemma
proved in~\cite{Count} (Lemma~\ref{counting2} below).
Full details can be found in~\cite{Olly}.
In order to state Lemma~\ref{counting2}, we need the following definition.
Given a complex~$\Hy$ with vertices~$x_1,\dots,x_t$, a
complex~$\G$ is called $(d_3,\delta_3,d_2,\delta_2,r,\Hy)$-regular if~$\G$
is $t$-partite with vertex classes~$V_1,\dots,V_t$ and satisfies the following properties:
\begin{itemize}
\item Let $P$ denote the underlying graph of~$\G$.
For every edge $x_ix_j\in E_2(\Hy)$ the bipartite graph $P[V_i,V_j]$
is $(d_2,\delta_2)$-regular.
\item For every hyperedge $e=x_hx_ix_j\in E_3(\Hy)$ there exists $d_e\ge d_3$
such that the triad~$P[V_h,V_i,V_j]$ is $(d_e,\delta_3,r)$-regular with respect
to the underlying hypergraph of~$\G$.
\end{itemize}
In this case, we say that a labelled copy~$H$ of~$\Hy$ in~$\G$ is
\emph{partition-respecting} if for all $i\in [t]$ the vertex of~$H$ corresponding
to~$x_i$ is contained in~$V_i$.

\begin{lemma}\label{counting2}
Let $t,r,n_0$ be positive integers and let
$\beta,d_2,d_3,\delta_2,\delta_3$ be positive constants such that
$$1/n_0\ll 1/r \ll\delta_2\ll\min\{\delta_3,d_2\}\le\delta_3\ll
\beta,d_3,1/t.$$ Then the following holds for all integers
$n\ge n_0$. Suppose that $\Hy$ is a complex with vertices~$x_1,\dots,x_t$.
Suppose also that~$\G$ is a $(d_3,\delta_3,d_2,\delta_2,r,\Hy)$-regular
complex with vertex classes $V_1,\dots,V_t$, all of size $n$.
Then~$\G$ contains at least $$(1-\beta)
n^{t} d_2^{e_2(\Hy)}\prod_{e\in E_3(\Hy)} d_e$$ labelled partition-respecting
copies of~$\Hy$.
\end{lemma}

Note that the difference to Lemma~\ref{countinglemma} is that Lemma~\ref{counting2}
only gives a lower bound and every vertex of $\Hy$ is to be embedded into a different
vertex class of $\G$. On the other hand, Lemma~\ref{counting2} allows for
different `hypergraph densities' between the clusters. (Actually, the proof below would 
permit this in Lemma~\ref{countinglemma} as well, see~\cite{Olly}.)%
     \COMMENT{Comment: It seems that apart from requiring that these densities are
all at least $d_3$ one also needs that they are at most $1-d_3$. Or alternatively,
one has to choose $\delta_3$ (and all the following constants) depending on all
these densities.}

To derive Lemma~\ref{countinglemma}, first assume that each of the vertex classes
$X_i$ of $\Hy$ contains exactly one vertex (i.e.~we want to embed every vertex of $\Hy$
into a different vertex class of $\G$).
In this case we only have to deduce the upper bound in Lemma~\ref{countinglemma}
from the lower bound in Lemma~\ref{counting2}.
As the (simple) proof of this is quite similar to the proof
for the complete case in~\cite{NR03}, we just describe the main idea here.
So consider any $\D \subseteq E_3(\Hy)$.
Now we construct a complex $\G_\D$ from $\G$ as follows. 
For any triple $hij \in \D$, we delete all hyperedges
from $\G[V_h,V_i,V_j]$ and add as hyperedges all those triangles contained
in the underlying graph induced by $V_h$, $V_i$ and $V_j$
which did not form a hyperedge in $\G[V_h,V_i,V_j]$.
(Thus $\G_\D$ may be viewed as a `partial' complement of $\G$.)
Now let $|\Hy^{(2)}|_\G$ denote the number of labelled partition-respecting
copies of the underlying graph of $\Hy$
in (the underlying graph of) $\G$. Then it is easy to see that 
$$
\sum_{\D \subseteq E(\Hy)} |\Hy|_{\G_\D} = |\Hy^{(2)}|_\G.
$$
(This is where we need to assume that we are considering the special case
when every vertex of $\Hy$ is embedded into a different vertex class of $\G$.)
We can use the (easy) counting lemma for graphs to estimate~$|\Hy^{(2)}|_\G$. 
Moreover, note that we are aiming for an upper bound on
the summand where $\D$ is empty. But we can obtain this since we can 
apply Lemma~\ref{counting2} 
to all the remaining summands. (This is where we need that Lemma~\ref{counting2}
allows for different `hypergraph densities'.)%
     \COMMENT{Comment: for this we need that $\delta_3\ll 1-d_3$. But the case when
$d_3=1$ is trivial. So $\delta_3\ll 1-d_3$ can be viewed to be implicit in
$\delta_3\ll 1-d_3$.}
A simple calculation gives the desired result.

So it remains to deduce the general case in Lemma~\ref{countinglemma}
from the special case when each of the vertex classes
$X_i$ of $\Hy$ contains exactly one vertex. 
To achieve this, consider the following construction.
Let $\G_1$ be the complex obtained from $\G$ by taking
$|X_1|$ copies of $\G$ and identifying them in $V(\G)\setminus V_1$. 
In other words, we blow up $V_1$ into $|X_1|$ copies, i.e.~$V_1$ 
is replaced by classes $V_{1i}$ with $1 \le i \le |X_1|$.
Now let $\G_2$ be the hypergraph obtained from $\G_1$ by taking
$|X_2|$ copies of $\G_1$ and identifying them in $V(\G_1)\setminus V_2$.
We continue in this way to obtain an $|\Hy|$-partite hypergraph $\G^*:=\G_k$
(so we have blown up each $V_i$ into $|X_i|$ copies).
Now view $\Hy$ as a $|\Hy|$-partite complex $\Hy^*$ with vertex classes
$X_{ji}$, each consisting of a single vertex, where $1 \le j \le k$ and $1 \le i \le |X_j|$.
Note that every labelled partition-respecting copy of~$\Hy$ in $\G$
yields a distinct labelled partition-respecting copy of~$\Hy^*$ in $\G^*$
(where in the latter case, $X_{ji}$ is mapped to $V_{ji}$).
So $|\Hy|_\G \le |\Hy^*|_{\G^*}$.
On the other hand, if a labelled partition-respecting copy of~$\Hy^*$ in $\G^*$
does not correspond to a labelled partition-respecting copy of~$\Hy$ in $\G$
then this means that this copy of $\Hy^*$ uses (at least) two `twin' vertices
in $\G^*$ which correspond to the same vertex in $\G$. There are at most
$|\Hy|n$ possibilities for choosing the first twin vertex, at most $|\Hy|$
possibilities for the second twin vertex and at most $n^{|\Hy|-2}$ possibilities
for the remaining vertices.
Thus  $|\Hy|_\G \ge |\Hy^*|_{\G^*}-|\Hy|^2 n^{|\Hy|-1}$.
We can now obtain the desired upper and lower bound on $|\Hy|_\G$
from the bounds on  $|\Hy^*|_{\G^*}$ which we already know.
(Note that these bounds imply that the `error term'
$|\Hy|^2 n^{|\Hy|-1}$ is negligible compared to  $|\Hy^*|_{\G^*}$.)


\section{Proof of Lemma~\ref{extensions_count}} \label{sec:extensions_count}

Throughout this section,
whenever we refer to copies of~$\Hy$ or~$\Hy'$ in~$\G$ we mean that these copies
will be labelled and partition-respecting without mentioning this explicitely.
As in Section~\ref{sec:pfemb} we denote such copies by~$H$ and~$H'$ respectively.
Thus, given any copy~$H$ of~$\Hy$, we have to estimate the number 
of extensions of~$H$ into copies of~$\Hy'$ in~$\G$. 
Recall that we denote the number of all these extensions by $|H \rightarrow \Hy'|_{\G}$.
Also, as in Section~\ref{sec:pfemb} we write
$$\overline{|\Hy\to\Hy'|}:=n^{t'-t} d_2^{e_2(\Hy')-e_2(\Hy)}d_3^{e_3(\Hy')-e_3(\Hy)}.
$$
We will argue 
similarly as in the proof of Corollary 26 in~\cite{roedlschacht}.
Namely, we use the following fact that can be deduced from the Cauchy-Schwartz inequality:
\begin{fact}\label{Cauchy-Schwartz}
For any $\beta>0$ there exists $\delta>0$, such that for any collection of non-negative
real numbers $x_1,\ldots , x_N$ satisfying
\begin{equation}
\sum_{i=1}^{N}x_i = (1\pm \delta)AN \text{\ \ \ and\ \ \ }
\sum_{i=1}^{N}x_i^2 = (1\pm \delta)A^2 N
\end{equation}
for some $A\geq 0$, all except at most $\beta N$ of the $x_i$'s lie in the interval 
$(1\pm \beta)A$.  
\end{fact}
\noindent
In our case, the collection $\{x_i\}_{i=1}^{N}$ will be $\{|H \rightarrow \Hy'|_{\G} \}_{H \in \G}$
(so $N:=|\Hy|_\G$) and we set
$$A:= \overline{|\Hy\to\Hy'|}.
$$
Given $\beta$ as in Lemma~\ref{extensions_count}, we let $\delta = \delta (\beta)$
be as in Fact~\ref{Cauchy-Schwartz}. We may assume that the hierarchy of constants in
Lemma~\ref{extensions_count} was chosen such that $\delta_3\ll\delta$.
To prove Lemma~\ref{extensions_count} it suffices to show that 
\begin{equation} \label{premise_one_too}
\sum_{H \in \G} |H \rightarrow \Hy'|_{\G} = (1\pm \delta)A |\Hy|_{\G},
\end{equation}
and 
\begin{equation} \label{premise_two} 
\sum_{H \in \G} |H \rightarrow \Hy'|_{\G}^2 = (1\pm \delta)A^2 |\Hy|_{\G} .
\end{equation}
The counting lemma (Lemma~\ref{countinglemma}) implies that 
\begin{equation*} 
|\Hy' |_{\G} = (1\pm \delta/8) n^{t'} d_2^{e_2(\Hy')}d_3^{e_3(\Hy')}
\end{equation*}
and
\begin{equation} \label{total_count_1} 
|\Hy |_{\G} = (1\pm \delta/8) n^{t} d_2^{e_2(\Hy)}d_3^{e_3(\Hy)}.
\end{equation}
It follows that
\begin{equation*} 
\sum_{H \in \G} |H \rightarrow \Hy'|_{\G} = |\Hy'|_{\G}=(1\pm \delta)A|\Hy|_{\G},
\end{equation*}
as required in~(\ref{premise_one_too}). 

To show~(\ref{premise_two}) we have to estimate
$\sum_{H \in \G} |H \rightarrow \Hy'|_{\G}^2$. Thus consider any
copy~$H$ of~$\Hy$ in~$\G$. Then $|H \rightarrow \Hy'|_{\G}^2$ corresponds to the number of
pairs $(H'_1,H'_2)$ of copies of $\Hy'$ in~$\G$ extending~$H$.
Now let~$\tilde{\Hy}'$ denote the complex which is
obtained from two disjoint copies of~$\Hy'$ by identifying them in~$V(\Hy)$.
Then the copies of~$\tilde{\Hy}'$ in~$\G$ which extend~$H$ correspond bijectively
to those pairs~$(H'_1,H'_2)$ which meet precisely in~$H$ and are disjoint otherwise.%
\footnote{Again, we only consider the partition-respecting copies
of~$\tilde{\Hy}'$ in~$\G$,
i.e.~if a vertex $\tilde{x}\in\tilde{\Hy}'$ corresponds to a vertex $x\in \Hy'$
which lies in~$X_i$, then~$\tilde{x}$ has to be embedded into~$V_i$.}
On the other hand, at most $(t'-t)^2n^{2(t'-t)-1}$ of the pairs $(H'_1,H'_2)$
meet in some vertex outside~$H$. Thus
\begin{equation}\label{eqsquare1} 
\sum_{H \in \G} |H \rightarrow \Hy'|_{\G}^2 \le
\sum_{H \in \G} \left(|H \rightarrow \tilde{\Hy}'|_{\G}+ (t'-t)^2n^{2(t'-t)-1}\right)
\leq |\tilde{\Hy}'|_{\G}+ (t'-t)^2n^{2t'-t-1}
\end{equation}
and clearly also
\begin{equation}\label{eqsquare2} 
\sum_{H \in \G} |H \rightarrow \Hy'|_{\G}^2\ge |\tilde{\Hy}'|_{\G}.
\end{equation}
But the counting lemma implies that
$$
|\tilde{\Hy}'|_{\G'} = (1\pm \delta/8) n^{2t'-t} d_2^{2e_2(\Hy')-e_2(\Hy)}d_3^{2e_3(\Hy')-e_3(\Hy)}
\stackrel{(\ref{total_count_1})}{=}(1\pm \delta/2)A^2|\Hy|_\G.
$$
In particular, $(t'-t)^2n^{2t'-t-1}\le \delta |\tilde{\Hy}'|_{\G'}/8\le \delta A^2|\Hy|_\G/2$.
Together with~(\ref{eqsquare1}) and~(\ref{eqsquare2}) this
implies~(\ref{premise_two}) and completes the proof of Lemma~\ref{extensions_count}.
Note that the proof above also allows for different `hypergraph densities' between the
clusters in Lemma~\ref{countinglemma}, but we have not included this to avoid making the
statement more technical.


\section{The regularity lemma for 3-uniform hypergraphs} \label{regularity}

\subsection{The Regularity Lemma -- definitions and statement}

The main purpose of this section is to introduce the regularity
lemma for 3-uniform hypergraphs due to Frankl and R\"odl~\cite{FR}.
As in the proof of the graph analogue of Theorem~\ref{Ramseythm}
we shall make use of it in order to obtain the necessary regular
complex~$\G$ to which we then apply the embedding lemma
(see Section~\ref{sec:proofofThm1} for the details). Before we can
state it, we will collect the necessary definitions.

\begin{definition}[$(\ell,t,\eps_1,\eps_2)$-partition]\label{defpartit}
{\rm Let $V$ be a set. An \emph{$(\ell,t,\eps_1,\eps_2)$-partition
$\Pa$ of $V$} is a partition into $V_0,V_1,\dots, V_t$ together with
families $(P_{\alpha}^{ij})_{\alpha=0}^{\ell_{ij}}$ ($1\le i<j\le
t$) of edge-disjoint bipartite graphs such that
\begin{itemize}
\item[{\rm (i)}] $|V_1|=\dots=|V_t|=\lfloor |V|/t\rfloor=:n$,
\item[{\rm (ii)}] $\ell_{ij}\le \ell$ for all pairs $1\le i<j\le t$,
\item[{\rm (iii)}] $\bigcup_{\alpha=0}^{\ell_{ij}} P_\alpha^{ij}$
is the complete bipartite graph with vertex classes $V_i$ and $V_j$
(for all pairs $1\le i<j\le t$),
\item[{\rm (iv)}] all but at most $\eps_1\binom{t}{2}n^2$ edges of the
complete $t$-partite graph $K[V_1,\dots,V_t]$ with vertex classes
$V_1,\dots,V_t$ lie in some $\eps_2$-regular graph $P_\alpha^{ij}$,
\item[{\rm (v)}] for all but at most $\eps_1\binom{t}{2}$ pairs $V_i,V_j$
($1\le i<j\le t$) we have $e(P_0^{ij})\le \eps_1 n^2$ and
$$| d_{P_{\alpha}^{ij}}(V_i,V_j)-1/\ell|\le \eps_2
$$
for all $\alpha=1,\dots,\ell_{ij}$.
\end{itemize}
}
\end{definition}

\begin{definition}[$(\delta_3,r)$-regular $(\ell,t,\eps_1,\eps_2)$-partition]\label{defregpartit}
{\rm Suppose that $\G$ is a $3$-uniform hypergraph and that
$V_0,V_1,\dots,V_t$ is an $(\ell,t,\eps_1,\eps_2)$-partition of the
vertex set $V(\G)$ of $\G$. Set $n:=|V_1|=\dots=|V_t|$. Recall that
a triad is a $3$-partite graph of the form $P=P_\alpha^{ij}\cup
P_\beta^{jk}\cup P_\gamma^{ik}$ and that $t(P)$ denotes the number of 
triangles in $P$. We say that the partition
$V_0,V_1,\dots,V_t$ is \emph{$(\delta_3,r)$-regular} if
$$ \sum_{\text{irregular}} t(P) < \delta_3|\G|^3,
$$
where $\sum_{\text{irregular}}$ denotes the sum over all triads $P$
which are not $(\delta_3,r)$-regular with respect to $\G$.}
\end{definition}

We can now state the regularity lemma for $3$-uniform hypergraphs
which was proved by Frankl and R\"odl~\cite{FR}.

\begin{theorem}[Regularity lemma for $3$-uniform hypergraphs]\label{3reglemma}
For all $\delta_3$ and $\eps_1$ with $0<\eps_1\le 2\delta_3^4$, for
all $t_0,\ell_0\in\mathbb{N}$ and for all integer-valued functions
$r=r(t,\ell)$ and all decreasing functions $\eps_2(\ell)$ with
$0<\eps_2(\ell)\le 1/\ell$, there exist integers $T_0, L_0$ and
$N_0$ such that the vertex set of any $3$-uniform hypergraph $\G$ of
order $|\G|\ge N_0$ admits a $(\delta_3,r)$-regular
$(\ell,t,\eps_1,\eps_2(\ell))$-partition for some $t$ and $\ell$
satisfying $t_0\le t\le T_0$ and $\ell_0\le\ell\le L_0$.
\end{theorem}

\noindent The elements $V_1,\dots,V_t$ of the
$(\ell,t,\eps_1,\eps_2(\ell))$-partition given by
Theorem~\ref{3reglemma} are called \emph{clusters}. $V_0$ is the
\emph{exceptional set}.

\subsection{Definition of the reduced hypergraph}\label{sec:rhgraph}

When we apply the graph regularity lemma to a graph $G$, we often
consider the so called reduced graph, whose vertices are the
clusters $V_i$ and whose edges correspond to those pairs of clusters
which induce
an $\eps$-regular bipartite graph. 
Analogously, we will now define a 3-uniform reduced hypergraph. 

In the proof of Theorem~\ref{Ramseythm} in the next section, we will
fix positive constants
satisfying the following hierarchy:
\begin{equation}\label{eqdefconst2}
\eps_1,1/t_0,1/\ell_0 \ll \delta_3\ll \eps_3 \ll 1/\Delta 
\end{equation}
where $\ell_0,t_0\in \mathbb{N}$ and we choose these constants
successively from right to left as explained earlier.
Next, for all $\ell\ge \ell_0$ and all $t\ge t_0$ we define
functions $r(t,\ell)$ and $\eps_2(\ell)$ satisfying the following
properties:
\begin{align} \label{eqdefeps2}
\frac{1}{r(t,\ell)} \ll \eps_2(\ell)\ll \frac{1}{\ell}, \delta_3,\eps_1.  
\end{align}
Suppose that with this choice of constants we have applied the regularity lemma
to a $3$-uniform hypergraph~$\G$. In particular, this gives an integer~$\ell$.
We then define constants~$d_2$ and~$\delta_2$ by
\begin{equation}\label{eqd2delta2}
d_2:=1/\ell \text{\ \ \ and \ \ \ } \delta_2:= \sqrt{\eps_2}.
\end{equation}
In order to define the reduced hypergraph corresponding to the partition
of~$V(\G)$ obtained from
the regularity lemma, we need the following definitions.
 
\begin{definition}[good pair $V_iV_j$]\label{defbadpair}
{\rm We call a pair $V_iV_j$ $(1\le i<j\le t$) of clusters
\emph{good} if it satisfies the following two properties:
\begin{itemize}
\item $e(P_0^{ij})\le \eps_1 n^2$ and
$| d_{P_{\alpha}^{ij}}(V_i,V_j)-d_2|\le \eps_2$ for all
$\alpha=1,\dots,\ell_{ij}$. (This means that $V_iV_j$ does not
belong to the at most $\eps_1\binom{t}{2}$ exceptional pairs
described in Definition~\ref{defpartit}(v).)
\item at most $\eps_3\ell/6$ of the bipartite graphs $P_\alpha^{ij}$
($1 \le \alpha\le \ell_{ij}$) are not
$(d_2,\delta_2)$-regular.
\end{itemize}
}
\end{definition}

\noindent 
Later on, we will use the fact that the first condition in Definition~\ref{defbadpair}
implies that $\ell_{ij}\ge \ell/2$ since~$d_2=1/\ell$.
An observation from~\cite{Ham} states that almost all pairs
of clusters are good, but we will not make use of this explicitely.

\begin{definition}[good triple $V_iV_jV_k$]\label{defbadtriple}
{\rm We call a triple $V_iV_jV_k$ $(1\le i<j<k\le t$) of clusters
\emph{good} if both of the following hold:
\begin{itemize}
\item each of the pairs $V_iV_j$, $V_jV_k$ and $V_iV_k$ is good,
\item at most $\eps_3\ell^3$ of the triads induced by $V_i,V_j,V_k$
are not $(\delta_3,r)$-regular with respect to~$\G$.
\end{itemize}
}
\end{definition}

The next proposition, which follows immediately from Proposition~5.12
in~\cite{Ham}, states that only a small
fraction of the triples $V_iV_jV_k$ are not good.

\begin{proposition}\label{badtriples}
At most $40\delta_3 \binom{t}{3}/\eps_3$ triples $V_iV_jV_k$ of clusters are
not good. 
\end{proposition}

We are now ready to define the reduced hypergraph $\R$.

\begin{definition}[Reduced hypergraph]\label{rhypergraph}
{\rm The vertices of the \emph{reduced hypergraph $\R$} are all the
clusters $V_1,\dots,V_t$. The hyperedges of $\R$ are precisely the
good triples $V_iV_jV_k$.}
\end{definition}

\noindent Thus, like $\G$, also $\R$ is a 3-uniform hypergraph.


\section{Proof of Theorem~\ref{Ramseythm}}\label{sec:proofofThm1}

In this section, we put together all the previous tools to prove
Theorem \ref{Ramseythm}. We will also make use of
the following well-known result (see e.g.~\cite{hyperturan}).

\begin{lemma}\label{turanlemma}
For all $k \in \mathbb{N}$ there exists a constant $c_0=c_0(k)<1$ such that if
$\mathcal{R}$ is a $3$-uniform hypergraph on $t \geq k$ vertices, and
if $e(\mathcal{R})\geq c_0 \binom{t}{3}$, then $\mathcal{R}$
contains a copy of~$K_k^{(3)}$.
\end{lemma}

We will also use the existence of hypergraph Ramsey numbers, without needing any explicit
upper bounds. Roughly speaking, the proof of Theorem~\ref{Ramseythm} proceeds
as follows. Consider any red/blue colouring of the hyperedges of $K_m^{(3)}$, where~$m$
is a sufficiently large integer (but~$m$ will be linear in~$|\Hy|$). We apply the
hypergraph regularity lemma to the red subhypergraph~$\G_{red}$ to obtain the reduced
hypergraph~$\mathcal{R}$, and show that~$\R$ satisfies the conditions of
Lemma~\ref{turanlemma} with $k:=R(K_{2\Delta+1}^{(3)})$.
Thus~$\mathcal{R}$ will contain a copy of~$K_k^{(3)}$. This copy corresponds
to~$k$ clusters such that for each triple of these clusters almost all the triads
are regular with respect to the red hypergraph~$\G_{red}$. We will then show that between each
pair $V_i,V_j$ of these clusters one can choose one of the bipartite graphs~$P^{ij}_\alpha$
in such a way that any triad $P_{hij}$ consisting of the chosen bipartite graphs is
regular with respect to~$\G_{red}$. Let~$P$ denote the~$k$-partite graph formed by all the
chosen bipartite graphs. We then consider the
following red/blue colouring of~$K_k^{(3)}$.
We colour the hyperedge~$hij$ with red if the triad $P_{hij}$ has density at least~$1/2$
with respect to~$\G_{red}$ and blue otherwise.
Since $k=R(K_{2\Delta+1}^{(3)})$ we can find a monochromatic
$K_{2\Delta+1}^{(3)}$. If it is red then we can apply the embedding lemma to 
the corresponding $(2\Delta+1)$-partite subhypergraph of~$\G_{red}$ and the corresponding
$(2\Delta+1)$-partite subgraph $P'$ of~$P$ to find a red
copy of~$\Hy$. This can be done since the chromatic number of~$\Hy$ is at
most~$2\Delta+1$ as~$\Delta(\Hy)\le \Delta$. If our monochromatic copy of
$K_{2\Delta+1}^{(3)}$ is blue then we can apply the embedding lemma to
the $(2\Delta+1)$-partite subhypergraph of the blue hypergraph~$\G_{blue} \subseteq K_m^{(3)}$ 
and~$P'$.

\bigskip

\noindent \textbf{Proof of Theorem \ref{Ramseythm}.} Let $m \in
\mathbb{N}$ be large enough for all subsequent calculations to hold.
We will check later that we can choose $m$ to be linear in $|\Hy|$.
Consider any red/blue-colouring of the hyperedges of $K_m^{(3)}$. 
Let~$\G_{red}$ be the red
and  $\G_{blue}$ be the blue subhypergraph on~$V(K_m^{(3)})$. We may assume
without loss of generality that $e(\G_{red})\geq e(\G_{blue})$.%
     \COMMENT{this is not necessary, but saves the reader from having to
check that it is ok to apply the RL to an extremly sparse graph}
We apply
the hypergraph regularity lemma to~$\G_{red}$ with parameters
$$t_0 \ge R(K_{2\Delta+1}^{(3)})=:k$$
as well as  $\ell_0,\delta_3,\eps_1$ and
functions $r(t,\ell)$ and $\eps_2(\ell)$ satisfying the hierarchies
in~(\ref{eqdefconst2}) and~(\ref{eqdefeps2}). 

Thus we obtain a set of clusters $V_1,\ldots,V_t$, each of size~$n$ say,
together with a partition~$(P^{ij}_\alpha)_{\alpha=0}^{\ell_{ij}}$
of the complete bipartite graph
between clusters~$V_i$ and~$V_j$ (for all $1\le i<j\le t$).
We define~$d_2$ and~$\delta_2$ as in~(\ref{eqd2delta2}) and let~$\R$
denote the reduced hypergraph. Proposition~\ref{badtriples} implies
that $\R$ has at least $(1-\eps)\binom{t}{3}$ hyperedges, where $\eps:=40\delta_3/\eps_3$.
Thus~(\ref{eqdefconst2}) implies
that $e(\mathcal{R}) \ge (1-\eps)\binom{|\mathcal{R}|}{3}> c_0 \binom{|\mathcal{R}|}{3}$,
where $c_0$ is as defined in Lemma~\ref{turanlemma}.
Since $|\mathcal{R}|\geq t_0 \ge k$, this means that we can apply Lemma~\ref{turanlemma}
to~$\mathcal{R}$ to obtain a copy of~$K_k^{(3)}$ in~$\mathcal{R}$.
Without loss of generality we may assume that the vertices of this copy are
the clusters $V_1,\dots,V_k$.

As indicated before, our next aim is to show that for each of the
$\binom{k}{2}$ pairs~$V_iV_j$ (with $1\le i<j \le k$)
one can choose one of the bipartite graphs
$P^{ij}_\alpha$ in such a way that each of them is $(d_2,\delta_2)$-regular
and such that each of the $\binom{k}{3}$ triads formed by the chosen
bipartite graphs is $(\delta_3,r)$-regular with respect to~$\G_{red}$.
We will denote the chosen bipartite graph between~$V_i$ and~$V_j$ by~$P_{ij}$ and
the triad between $V_h$, $V_i$ and $V_j$ by~$P_{hij}$.

To see that such graphs $P_{ij}$ exist, consider selecting (for each pair $i,j$)
one of the $\ell_{ij}$ bipartite graphs $P_\alpha^{ij}$ with
$1\le \alpha \le \ell_{ij}$ uniformly at random.
By Definition~\ref{defbadpair}, the probability that $P_{ij}$
is not $(d_2,\delta_2)$-regular is at most $(\eps_3 \ell/6)/\ell_{ij} \le \eps_3/3$.
So the probability that all of the selected bipartite graphs are
$(d_2,\delta_2)$-regular is at least 
\begin{equation} \label{probbound}
1- \binom{k}{2}\frac{\eps_3}{3} \stackrel{(\ref{eqdefconst2})}{>} \frac{1}{2}.
\end{equation} 
Similarly, Definition~\ref{defbadtriple} implies that the probability that the 
triad $P_{hij}$ is not $(\delta_3,r)$-regular is at most
$\eps_3 \ell^3/\ell_{hi} \ell_{ij} \ell_{hj} \le 8\eps_3$.
So the probability that all of the selected $P_{hij}$ are $(\delta_3,r)$-regular
is at least $1-\binom{k}{3}8\eps_3 > 1/2$.
Together with~(\ref{probbound}), this shows that there is some choice of 
bipartite graphs $P_{ij}$ which has the required properties.

We now use the densities of the corresponding triads $P_{hij}$ to
define a red/blue-colouring of the $K_k^{(3)}$ which we found in $\R$:
if $d_{\G_{red}} (P_{hij}) \ge 1/2$, then we colour the hyperedge $V_hV_iV_j$ red, 
otherwise we colour it blue. 
Since $k=R(K_{2\Delta+1}^{(3)})$, we find a monochromatic copy $K$ of 
$K_{2\Delta+1}^{(3)}$ in our $K_k^{(3)}$. We
now greedily assign the vertices of $\Hy$ to the clusters that form
the vertex set of $K$, in such a way that if three vertices of
$\Hy$ form a hyperedge, then they are assigned to different
clusters. (We may think of this as a $(2\Delta +1)$-vertex-colouring
of $\Hy$.)
We now need to show that with this assignment we can apply the
embedding lemma to find a monochromatic copy of $\Hy$ in~$K_m^{(3)}$.

Assume first that $K$ is red. We already have
bipartite graphs~$P_{ij}$ between the clusters in $V(K)$ which are
$(d_2,\delta_2)$-regular and form triads~$P_{hij}$ which are
$(\delta_3,r)$-regular with respect to~$\G_{red}$. The only technical
problem is
that these triads do not all have the same density with respect to $\G_{red}$,
which was one of the conditions in the embedding lemma. We do know,
however, that in each case we have $d_{\G_{red}}(P_{hij}) \geq 1/2$. 
So we choose a hypergraph $\G_{red}'\subseteq \G_{red}$
such that all the graph triads are $(1/2,3\delta_3,r)$-regular with
respect to $\G_{red}'$.
It is easy to see that such a $\G_{red}'$ exists:
for each triple $V_hV_iV_j$ that is a hyperedge of $K$, 
consider a random subset of the hyperedges of $\G_{red}$
induced by $V_i,V_j,V_k$ such that $P_{hij}$ has density $(1\pm \delta_3)/2$
with respect to this subset.
This observation is formalized for instance 
in Proposition~33 of~\cite{roedlschacht}, which%
    \COMMENT{Comment: It's a bit unfortunate that the proof is
actually not in that paper, but since the result really is quite easy,
it's not worth giving further details here}
one can apply directly to obtain the
above bounds on the regularity of $\G_{red}'$.
(Alternatively, it is easy to see that the proof of
Lemma~\ref{emblemma} generalizes to different `hypergraph densities'.)
We then apply the embedding lemma (Lemma~\ref{embcor}) to find
a copy of $\Hy$ in $\G_{red}'$, and therefore also in $\G_{red}$.

On the other hand, if $K$ is blue, we will aim to
find a copy of $\Hy$ in $\G_{blue}$. We certainly still have a set of
bipartite graphs all of which are $(d_2,\delta_2)$-regular, but we
now also need to prove that all triads are regular with respect to
$\G_{blue}$.
So suppose $\vec{Q}=(Q(1),\ldots,Q(r))$ is an
$r$-tuple of subtriads of one of these triads~$P_{hij}$, satisfying
$t(\vec{Q})>\delta_3t(P_{hij})$. Let~$d$ be such that~$P_{hij}$ is
$(d,\delta_3,r)$-regular with respect to $\G_{red}$. Then
\begin{align*}|(1-d)-d_{\G_{blue}}(\vec{Q})|
 = |d-(1-d_{\G_{blue}}(\vec{Q}))|
 = |d-d_{\G_{red}}(\vec{Q})|
 < \delta_3.
\end{align*}
Thus $P_{hij}$ is $(1-d,\delta_3,r)$-regular with respect to~$\G_{blue}$
(note that $\delta_3 \ll  1/2 \le 1-d$). By
the same method as in the previous case, we can apply the embedding
lemma to obtain a copy of $\Hy$ in $\G_{blue}$.

It remains to estimate how large we needed $m$ to be in order
for all of our calculations to be valid. When we apply the embedding
lemma, we know we can find any subgraph $\Hy$ of maximum degree at most~$\Delta$
with $|\Hy|\le cn$, where~$n$ is the size of a cluster and~$c$ is a
constant chosen to satisfy the conditions of the embedding lemma.
Since $n=\lfloor m/t\rfloor\ge m/2T_0$, this means that it suffices
to start with an $m$ satisfying $m\ge 2T_0|\Hy|/c$.
In order to be able to apply the embedding lemma we need that
$c \ll d_2,d_3=1/2,1/\Delta$. We obtain $d_2=1/\ell$ from the regularity lemma,
given constants $\delta_3,\eps_1,t_0,\ell_0$, an integer-valued
function $r=r(t,\ell)$ and a decreasing function $\eps_2(\ell)$,
all satisfying the hierarchies~(\ref{eqdefconst2}) and~(\ref{eqdefeps2}).
In all cases, we can view the constants we require purely as functions of $\Delta$.
Thus~$c$ is implicitly a function solely of $\Delta$.
This is also the case for~$T_0$ and~$N_0$.

Finally, in order to be able to apply the regularity lemma to $\G_{red}$ we
needed to assume that $m \geq N_0$, and in order to be able
to apply the embedding lemma we needed to assume that $n \geq n_0$
(for which it is sufficient to assume that $m \geq 2T_0n_0$).
Altogether, this
shows that we can take the constant~$C$ in Theorem~\ref{Ramseythm}
to be $\max \{2T_0/c,N_0,2T_0n_0 \}$.
\noproof


\medskip

{\footnotesize \obeylines \parindent=0pt

Oliver Cooley, Nikolaos Fountoulakis, Daniela K\"{u}hn \& Deryk Osthus 
School of Mathematics
University of Birmingham
Edgbaston
Birmingham
B15 2TT
UK
}

{\footnotesize \parindent=0pt

\it{E-mail addresses}:
\tt{\{cooleyo,nikolaos,kuehn,osthus\}@maths.bham.ac.uk}}

\end{document}